  \documentstyle[12pt]{article}
\newcommand{\numero}[1]{
\addtocounter{section}{1}
\begin{center}{\bf \thesection .\
#1\vspace{-.1in}}\end{center}
\setcounter{subsection}{0}
\setcounter{lemma}{0}\indent}

\newcommand{\subnumero}[1]{
\pagebreak[1]\begin{center}{\em #1}\nopagebreak\end{center}
}

\newtheorem{lemma}{Lemma}[section]
\newtheorem{theorem}[lemma]{Theorem}
\newtheorem{definition}[lemma]{Definition}
\newtheorem{corollary}[lemma]{Corollary}

\newtheorem{conjecture}{Conjecture}
\newtheorem{proposition}[lemma]{Proposition}

\newcommand{\rr}{{\bf R}}

\newcommand{\zz}{{\bf Z}}

\newcommand{\eop}{\hfill $/$\hspace*{-.1cm}$/$\hspace*{-.1cm}$/$\vspace{.1in}}

\begin{document}

\section*{Homotopy types of strict $3$-groupoids}

Carlos Simpson\newline
CNRS, UMR 5580, Universit\'e de Toulouse 3

\bigskip

It has been difficult to see precisely the role played by {\em strict}
$n$-categories in the nascent theory of $n$-categories, particularly as related
to $n$-truncated homotopy types of spaces. We propose to show in a fairly
general setting that one cannot obtain all $3$-types by any reasonable
realization functor
\footnote{
Our notion of ``reasonable realization functor'' (Definition
\ref{realizationdef}) is any functor $\Re$ from the category of strict
$n$-groupoids to $Top$, provided with a natural transformation $r$ from the set
of objects of $G$ to the points of $\Re (G)$, and natural isomorphisms $\pi
_0(G)\cong \pi _0(\Re (G))$ and $\pi _i(G,x) \cong \pi _i(\Re (G), r(x))$.
This axiom is fundamental to the question of whether one can realize homotopy
types by strict $n$-groupoids, because one wants to read off the homotopy
groups of the space from the strict $n$-groupoid. The standard
realization functors satisfy this property, and the somewhat different
realization construction of \cite{KV} is claimed there to have this property.
}
from strict $3$-groupoids (i.e. groupoids in the sense of
\cite{KV}).  More precisely we show that one does not obtain
the $3$-type of $S^2$. The basic reason is that the Whitehead bracket is
nonzero. This phenomenon is actually well-known, but in order to take into
account the possibility of an arbitrary reasonable realization functor we have
to write the argument in a particular way.

We start by recalling the notion of strict $n$-category. Then we look at the
notion of strict $n$-groupoid as defined by Kapranov and Voevodsky \cite{KV}.
We show that their definition is equivalent to a couple of other
natural-looking definitions (one of these equivalences was left as an exercise
in \cite{KV}). At the end of these first sections, we have a picture of strict
$3$-groupoids having only one object and one $1$-morphism, as being equivalent
to abelian monoidal objects $(G,+)$ in the category of groupoids, such that
$(\pi _0(G),+)$ is a group. In the case in question, this group will be $\pi
_2(S^2)=\zz$. Then comes the main part of the argument. We show that, up to
inverting a few equivalences, such an object has a morphism giving a splitting
of the Postnikov tower (Proposition \ref{diagramme}. It follows that for any
realization functor respecting homotopy groups, the Postnikov tower of the
realization (which has two stages corresponding to $\pi _2$ and $\pi _3$)
splits. This implies that the $3$-type of $S^2$ cannot occur as a realization.

The fact that strict $n$-groupoids are not appropriate for modelling all
homotopy types has in principle been known for some time. There are several
papers by R. Brown and coauthors on this subject, see \cite{RBrown1},
\cite{BrownGilbert}, \cite{BrownHiggins}, \cite{BrownHiggins2}; a recent paper
by C. Berger \cite{Berger}; and also a discussion of this in various places in
Grothendieck \cite{Grothendieck}. Other related examples are given in
Gordon-Power-Street \cite{Gordon-Power-Street}.
The novelty of our present treatment is that we have written the argument
in such a way that it applies to a wide class of possible realization
functors, and in particular it applies to the realization functor of
Kapranov-Voevodsky (1991) \cite{KV}.

This problem with strict $n$-groupoids can be summed up by saying in R.
Brown's terminology, that they correspond to {\em crossed complexes}. While
a nontrivial  action of $\pi _1$ on the $\pi _i$ can occur in a crossed complex,
the higher Whitehead operations such as $\pi _2\otimes \pi _2\rightarrow
\pi _3$
must vanish.  This in turn is due to the fundamental ``interchange rule'' (or
``Godement relation'' or ``Eckmann-Hilton argument''). This effect occurs when
one takes two $2$-morphisms $a$ and $b$ both with source and target a
$1$-identity  $1_x$. There are various ways of composing $a$ and $b$ in this
situation, and comparison of these compositions leads to the conclusion that all
of the compositions are commutative. In a weak $n$-category, this commutativity
would only hold up to higher homotopy, which leads to the notion of
``braiding''; and in fact it is exactly the braiding which leads to the
Whitehead operation. However, in a strict $n$-category, the commutativity is
exact, so the Whitehead operation is trivial.

One can observe that one of the reasons why this problem occurs is that we have
the exact $1$-identity $1_x$. This leads to wondering if one could get a better
theory by getting rid of the exact identities. We speculate in this direction
at the end of the paper by proposing a notion of {\em $n$-snucategory}, which
would be an $n$-category with strictly associative composition, but without
units; we would only require existence of weak units. The details of the notion
of weak unit are not worked out.

A preliminary version of this note was circulated in a limited way in the summer
of 1997.

I would like to thank: R. Brown, A. Brugui\`eres, A. Hirschowitz, G.
Maltsiniotis, and Z. Tamsamani.

\numero{Strict $n$-categories}

In what follows {\em all $n$-categories are meant
to be strict $n$-categories}. For this reason we try to put in the adjective
``strict'' as much as possible when $n>1$; but in any case, the very few
times that we speak of weak $n$-categories, this will be explicitly stated.
We mostly restrict our attention to $n\leq 3$.

In case that isn't already clear, it should be stressed that everything we do
in this section (as well as most of the next and even the subsequent one as
well) is very well known and classical, so much so that I don't know what are
the original references.

To start with, a {\em
strict $2$-category} $A$ is a collection of objects $A_0$ plus, for each pair of
objects $x,y\in A_0$ a category $Hom _A(x,y)$ together with a morphism
$$
Hom_A(x,y)\times Hom _A(y,z)\rightarrow Hom _A(x,z)
$$
which is strictly associative in the obvious way; and such that a unit exists,
that is an element $1_x\in Ob\, Hom _A(x,x)$ with the property that
multiplication by $1_x$ acts trivially on objects of $Hom _A(x,y)$ or
$Hom_A(y,x)$ and multiplication by $1_{1_x}$ acts trivially on morphisms of
these categories.

A {\em strict $3$-category} $C$ is the same as above but where $Hom _C(x,y)$ are
supposed to be strict $2$-categories. There is an obvious notion of direct
product of
strict $2$-categories, so the above definition applies {\em mutatis mutandis}.

For general $n$, the well-known definition is most
easily presented by induction on $n$. We assume known the definition of strict
$n-1$-category for $n-1$, and we assume known that the category of strict
$n-1$-categories is closed under direct product. A {\em strict $n$-category}
$C$ is then a category enriched \cite{Kelly} over the category of strict
$n-1$-categories. This means that $C$ is composed of a {\em set of objects}
$Ob(C)$ together with, for each pair $x,y\in Ob(C)$, a {\em morphism-object}
$Hom _C(x,y)$ which is a strict $n-1$-category; together with a strictly
associative composition law
$$
Hom _C(x,y)\times Hom _C(y,z) \rightarrow Hom _C(x,z)
$$
and a morphism $1_x: \ast \rightarrow Hom _C(x,x)$ (where $\ast$ denotes the
final object cf below) acting as the identity for the composition law. The {\em
category of strict $n$-categories} denoted $nStrCat$ is the category whose
objects are as above and whose morphisms are the transformations strictly
perserving all of the structures. Note that $nStrCat$ admits a direct product:
if $C$ and $C'$ are two strict $n$-categories then $C\times C'$ is the strict
$n$-category with
$$
Ob(C\times C'):= Ob(C) \times Ob(C')
$$
and for $(x,x'), \; (y,y') \in Ob(C\times C')$,
$$
Hom _{C\times C'}((x,x'), (y,y')):= Hom _C(x,y)\times Hom _{C'}(x',y')
$$
where the direct product on the right is that of $(n-1)StrCat$. Note that the
final object of $nStrCat$ is the strict $n$-category $\ast$ with exactly one
object $x$ and with $Hom _{\ast}(x,x)= \ast$ being the final object of
$(n-1)StrCat$.

The induction inherent in this definition may be worked out explicitly to give
the definition as it is presented in \cite{KV} for example. In doing this one
finds that underlying a strict $n$-category $C$ are the sets $Mor ^i(C)$ of {\em
$i$-morphisms} or {\em $i$-arrows}, for $0\leq i\leq n$. The $0$-morphisms
are by
definition the objects, and $Mor ^i(C)$ is the disjoint union over all
pairs $x,y$ of the $Mor ^{i-1}(Hom
_C(x,y))$. The composition laws at each stage lead to various compositions
for $i$-morphisms, denoted in \cite{KV} by $\ast _j$ for $0\leq j < i$.
These are partially defined depending on the {\em source} and {\em target} maps.
For a more detailed explanation, refer to the standard references
\cite{BrownHiggins} \cite{Street} \cite{KV} (and I am probably missing many
older references which could date back even before \cite{Benabou}
\cite{GabrielZisman}).

One of the most important of the axioms satisfied by the various compositions
in a strict $n$-category is variously known under the name of ``Eckmann-Hilton
argument'', ``Godement relations'', ``interchange rules'' etc.
The following discussion of this axiom owes a lot to discussions I had with
Z. Tamsamani during his thesis work. This axiom  comes from the fact that the
composition law
$$
Hom _C(x,y)\times Hom _C(y,z)\rightarrow Hom _C(x,z)
$$
is a morphism with domain the direct product of the two morphism
$n-1$-categories from $x$ to $y$ and from $y$ to $z$. In a direct product,
compositions in the two factors by definition are independent (commute).
Thus, for $1$-morphisms in $Hom _C(x,y)\times Hom _C(y,z)$ (where the
composition $\ast _0$ for these $n-1$-categories is actually the composition
$\ast _1$ for $C$ and we adopt the latter notation),
we have
$$
(a,b) \ast _1 (c,d) = (a\ast _1c, b \ast _1d).
$$
This leads to the formula
$$
(a\ast _0b) \ast _1 (c\ast _0d) = (a\ast _1c) \ast _0 (b \ast _1d).
$$
This seemingly innocuous formula takes on a special meaning when we start
inserting identity maps. Suppose $x=y=z$ and let $1_x$ be the identity of $x$
which may be thought of as an object of $Hom _C(x,x)$. Let $e$ denote the
$2$-morphism of $C$, identity of $1_x$; which may be thought of as a
$1$-morphism of $Hom _C(x,x)$. It acts as the identity for both compositions
$\ast _0$ and $\ast _1$ (the reader may check that this follows from the part
of the axioms for an $n$-category saying that the morphism $1_x: \ast
\rightarrow Hom _C(x,x)$ is an identity for the composition).

If $a, b$ are also
endomorphisms of $1_x$, then the above rule specializes to:
$$
a\ast _1b =
(a\ast
_0e) \ast _1 (e\ast _0b) = (a\ast _1e) \ast _0 (e \ast _1b)  = a \ast _0b.
$$
Thus in this case the compositions $\ast _0$ and $\ast _1$ are the same.
A different ordering gives the formula
$$
a\ast _1b =
(e\ast
_0a) \ast _1 (b\ast _0e) = (e\ast _1b) \ast _0 (a \ast _1e)  = b \ast _0a.
$$
Therefore we have
$$
a\ast _1b= b\ast _1a = a \ast _0b = b\ast _0a.
$$
This argument says, then, that $Ob (Hom _{Hom _C(x,x)}(1_x, 1_x))$
is a commutative monoid and the two natural multiplications are the same.

The same argument extends to the whole monoid structure on the $n-2$-category
$Hom _{Hom _C(x,x)}(1_x, 1_x)$:

\begin{lemma}
\label{godement}
The two composition laws on the strict $n-2$-category
$Hom _{Hom _C(x,x)}(1_x, 1_x)$ are equal, and this law is commutative.
In other words, $Hom _{Hom _C(x,x)}(1_x, 1_x)$ is an abelian monoid-object in
the category $(n-2)StrCat$.
\end{lemma}
\eop

There is a partial converse to the above observation: if the only object is $x$
and the only $1$-morphism is $1_x$ then nothing else can happen and we get the
following equivalence of categories.

\begin{lemma}
\label{scholium}
Suppose $G$ is an abelian
monoid-object in the category $(n-2)StrCat$. Then there is a unique strict
$n$-category $C$ such that
$$
Ob(C)= \{ x\}\;\;\; \mbox{and}\;\;\; Mor ^1(C)=Ob(Hom _C(x,x))=\{
1_x\}
$$
and such that $Hom _{Hom _C(x,x)}(1_x, 1_x)=G$ as an abelian
monoid-object. This construction establishes an equivalence between the
categories of abelian monoid-objects in $(n-2)StrCat$, and the strict
$n$-categories having only one object and one $1$-morphism.
\end{lemma}
{\em Proof:}
Define the strict $n-1$-category $U$ with $Ob(U)= \{ u\}$ and $Hom _U(u,u)=G$
with its monoid structure as composition law. The fact that the composition law
is commutative allows it to be used to define an associative and commutative
multiplication
$$
U\times U \rightarrow U.
$$
Now let $C$ be the strict $n$-category
with $Ob(C)=\{ x\}$ and $Hom _C(x,x)=U$ with the above multiplication. It is
clear that this construction is inverse to the previous one.
\eop

It is clear from the construction (the fact that the multiplication on $U$ is
again commutative) that the construction can be iterated any number of times.
We obtain the following corollary.

\begin{corollary}
\label{iterate}
Suppose $C$ is a strict $n$-category with only one object and only one
$1$-morphism. Then there exists a strict $n+1$-category $B$ with only one
object $b$ and with $Hom _B(b,b)\cong C$.
\end{corollary}
{\em Proof:}
By the previous lemmas, $C$ corresponds to an abelian monoid-object $G$ in
$(n-2)StrCat$. Construct $U$ as in the proof of \ref{scholium}, and note that
$U$ is an abelian monoid-object in $(n-1)StrCat$. Now apply the result of
\ref{scholium} directly to $U$ to obtain $B\in (n+1)StrCat$, which will have
the desired property.
\eop

\numero{The groupoid condition}

Recall that a {\em groupoid} is a category where all morphisms are invertible.
This definition generalizes to strict $n$-categories in the following way
\cite{KV}. We give a theorem stating that three versions of this definition are
equivalent.

Note that, following \cite{KV}, we {\em do not} require strict invertibility of
morphisms, thus the notion of strict $n$-groupoid is more general than the
notion employed by Brown and Higgins \cite{BrownHiggins}.

Our discussion is in many ways parallel to the treatment of the groupoid
condition for weak $n$-categories in \cite{Tamsamani} and our treatment in this
section comes in large part from discussions with Z. Tamsamani about this.

The statement of the theorem-definition is recursive on $n$.

\begin{theorem}
\label{thmdef}
Fix $n<\infty$.

\noindent
{\bf I. Groupoids}\,\,
Suppose $A$ is a strict $n$-category. The following three conditions
are equivalent (and in this case we say that $A$ is a {\em strict
$n$-groupoid}). \newline
(1)\, $A$ is an $n$-groupoid in the sense of Kapranov-Voevodsky
\cite{KV};
\newline
(2)\, for all $x,y\in A$, $Hom _A(x,y)$ is a strict $n-1$-groupoid, and for any
$1$-morphism $f:x\rightarrow y$ in $A$, the two morphisms of composition with
$f$
$$
Hom _A(y,z)\rightarrow Hom _A(x,z),\;\;\;\;
Hom _A(w,x)\rightarrow Hom_A(w,y)
$$
are equivalences of strict $n-1$-groupoids (see below);
\newline
(3)\, for all $x,y\in A$, $Hom _A(x,y)$ is a strict $n-1$-groupoid, and
$\tau _{\leq
1}A$ (defined below) is a $1$-groupoid.

\noindent
{\bf II. Truncation}\,\,
If $A$ is a strict $n$-groupoid, then define $\tau _{\leq k}A$ to be the
strict $k$-category whose $i$-morphisms are those of $A$ for $i<k$ and whose
$k$-morphisms are the equivalence classes of $k$-morphisms of $A$ under the
equivalence relation that two are equivalent if there is a $k+1$-morphism
joining them. The fact that this is an equivalence relation is a statement about
$n-k$-groupoids. The set $\tau _{\leq 0}A$ will also be denoted $\pi _0A$.
The truncation is again a $k$-groupoid, and for $n$-groupoids $A$ the truncation
coincides with the operation defined in \cite{KV}.

\noindent
{\bf III. Equivalence}\,\,
A morphism $f:A\rightarrow B$ of strict $n$-groupoids is said to be an {\em
equivalence} if the following equivalent conditions are satisfied:
\newline
(a)\, (this is the definition in
\cite{KV}) $f$ induces an isomorphism $\pi _0A\rightarrow
\pi _0B$, and for every object $a\in A$ $f$ induces an isomorphism
$\pi _i(A,a)\stackrel{\cong}{\rightarrow} \pi _i(B, f(a))$ where these homotopy
groups are as defined in \cite{KV};
\newline
(b)\, $f$ induces a surjection $\pi _0A\rightarrow \pi _0B$ and for every pair
of objects $x,y\in A$ $f$ induces an equivalence of $n-1$-groupoids
$Hom_A(x,y)\rightarrow Hom _B(f(x), f(y))$;
\newline
(c)\, if $u,v$ are $i$-morphisms in $A$ sharing the same source and target, and
if $r$ is an $i+1$-morphism in $B$ going from $f(u)$ to $f(v)$ then there exists
an $i+1$-morphism $t$ in $A$ going from $u$ to $v$ and an $i+2$-morphism in $B$
going from $f(t)$ to $r$ (this includes the limiting cases $i=-1$ where $u$ and
$v$ are not specified, and $i=n-1, n$ where ``$n+1$-morphisms'' mean equalities
between $n$-morphisms and ``$n+2$-morphisms'' are not specified).

\noindent
{\bf IV. Sub-lemma}\,\, If $f: A\rightarrow B$ and $g: B\rightarrow C$ are
morphisms of strict $n$-groupoids and if any two of $f$, $g$ and $gf$ are
equivalences, then so is the third.

\noindent
{\bf V. Second sub-lemma}\,\, If
$$
A\stackrel{f}{\rightarrow}B\stackrel{g}{\rightarrow}C\stackrel{h}{\rightarrow}D
$$
are morphisms of strict $n$-groupoids and if $hg$ and $gf$ are equivalences,
then $g$ is an equivalence.
\end{theorem}
{\em Proof:} It is clear for $n=0$, so we assume $n\geq 1$ and
proceed by induction on $n$: we assume that the theorem is true (and all
definitions are known) for  strict $n-1$-categories.

We first discuss the existence of truncation (part II), for $k\geq 1$. Note that
in this case $\tau _{\leq k}A$ may be defined as the strict $k$-category
with the
same objects as $A$ and with
$$
Hom _{\tau _{\leq k}A}(x,y):= \tau _{\leq k-1}Hom _A(x,y).
$$
Thus the fact that the relation in question is an equivalence relation, is a
statement about $n-1$-categories and known by induction. Note that the
truncation operation clearly preserves any one of the three groupoid conditions
(1), (2), (3). Thus we may affirm in a strong sense that
$\tau _{\leq k}(A)$ is a
$k$-groupoid without knowing the equivalence of the
conditions (1)-(3).

Note also that the truncation operation for $n$-groupoids is the same as that
defined in \cite{KV} (they define truncation for general strict
$n$-categories but for $n$-categories which are not groupoids, their definition
is different from that of \cite{Tamsamani} and not all that useful).

For $0\leq k\leq k'\leq n$ we have
$$
\tau _{\leq k}(\tau _{\leq k'}(A)) = \tau _{\leq k}(A).
$$
To see this note that the equivalence relation used to define the $k$-arrows of
$\tau _{\leq k}(A)$ is the same if taken in $A$ or in $\tau _{\leq
k+1}(A)$---the existence of a $k+1$-arrow going between two $k$-arrows
is equivalent to the existence of an equivalence class of $k+1$-arrows going
between the two $k$-arrows.

Finally using the above remark we obtain the existence of the truncation $\tau
_{\leq 0}(A)$: the relation is the same as for the truncation $\tau _{\leq
0}(\tau _{\leq 1}(A))$, and $\tau _{\leq 1}(A)$ is a strict $1$-groupoid in
the usual sense so the arrows are invertible, which shows that the relation
used to define the $0$-arrows (i.e. objects) in $\tau _{\leq 0}(A)$ is
in fact an equivalence relation.

We complete our discussion of truncation
by noting that there is a natural morphism
of strict $n$-categories
$A\rightarrow \tau _{\leq k}(A)$, where the right hand side ({\em a priori} a
strict $k$-category) is considered as a strict $n$-category in the obvious way.

\bigskip

We turn next to the notion of equivalence (part III), and prove that conditions
(a) and (b) are equivalent. This notion for $n$-groupoids will not enter
into the
subsequent treatment of part (I)---what does enter is the notion of equivalence
for $n-1$-groupoids, which is known by induction---so we may assume the
equivalence of definitions (1)-(3) for our discussion of part (III).

Recall first of all the definition of the homotopy groups. Let $1^i_a$ denote
the $i$-fold iterated identity of an object $a$; it is an $i$-morphism, the
identity of $1^{i-1}_a$ (starting with $1^0_a=a$). Then
$$
\pi _i(A,a):= Hom _{\tau _{\leq i}(A)}(1^{i-1}_a, 1^{i-1}_a).
$$
This definition is completed by setting $\pi _0(A):= \tau _{\leq 0}(A)$.
These definitions are the same as in \cite{KV}. Note directly from the
definition that for $i\leq k$  the truncation morphism induces isomorphisms
$$
\pi _i(A,a)\stackrel{\cong}{\rightarrow}\pi _i(\tau _{\leq k}(A), a).
$$
Also for $i\geq 1$ we have
$$
\pi _i(A,a)= \pi _{i-1}(Hom _A(a,a), 1_a).
$$
One shows that the $\pi _i$ are abelian for $i\geq 2$.  This is part of a more
general principle, the
``interchange rule'' or ``Godement relations'' refered to in \S 1.

Suppose $f:A\rightarrow B$ is a morphism of strict $n$-groupoids satisfying
condition (b). From the immediately preceding formula and the inductive
statement for $n-1$-groupoids, we get that $f$ induces isomorphisms
on the $\pi _i$ for $i\geq 1$. On the other hand, the truncation
$\tau _{\leq 1}(f)$ satisfies condition (b) for a morphism of $1$-groupoids,
and this is readily seen to imply that $\pi _0(f)$ is an isomorphism. Thus $f$
satisfies condition (a).

Suppose on the other hand that $f:A\rightarrow B$
is a morphism of strict $n$-groupoids satisfying
condition (a). Then of course $\pi _0(f)$ is surjective. Consider two objects
$x,y\in A$ and look at the induced morphism
$$
f^{x,y}: Hom _A(x,y)\rightarrow Hom _B(f(x), f(y)).
$$
We claim that $f^{x,y}$ satisfies condition (a) for a morphism of
$n-1$-groupoids. For this, consider a $1$-morphism from $x$ to $y$, i.e. an
object $r\in Hom _A(x,y)$. By version (2) of the groupoid condition
for $A$, multiplication by $r$ induces an equivalence of $n-1$-groupoids
$$
m(r): Hom _A(x,x)\rightarrow Hom _A(x,y),
$$
and furthermore $m(r)(1_x)=r$. The same is true in $B$: multiplication by
$f(r)$ induces an equivalence
$$
m(f(r)): Hom _B(f(x), f(x))\rightarrow Hom _B(f(x), f(y)).
$$
The fact that $f$ is a morphism implies that these fit into a commutative
square
$$
\begin{array}{ccc}
Hom _A(x,x)&\rightarrow &Hom _A(x,y)\\
\downarrow && \downarrow \\
Hom _B(f(x), f(x))&\rightarrow &Hom _B(f(x), f(y)).
\end{array}
$$
The equivalence condition (a) for $f$ implies that the left vertical morphism
induces isomorphisms
$$
\pi _i(Hom _A(x,x), 1_x)\stackrel{\cong}{\rightarrow}
\pi _i(Hom _B(f(x), f(x)), 1_{f(x)}).
$$
Therefore the right vertical morphism (i.e. $f_{x,y}$) induces isomorphisms
$$
\pi _i(Hom _A(x,y), r)\stackrel{\cong}{\rightarrow}
\pi _i(Hom _B(f(x), f(y)), f(r)),
$$
this for all $i\geq 1$.
We have now verified these isomorphisms for any base-object $r$. A similar
argument implies that $f^{x,y}$ induces an injection on $\pi _0$. On the other
hand, the fact that $f$ induces an isomorphism on $\pi _0$ implies that
$f^{x,y}$ induces a surjection on $\pi _0$ (note that these last two statements
are reduced to statements about $1$-groupoids by applying $\tau _{\leq 1}$ so we
don't give further details). All of these statements taken together  imply that
$f^{x,y}$ satisfies condition (a), and by the inductive statement of
the theorem for $n-1$-groupoids this implies that $f^{x,y}$ is an equivalence.
Thus $f$ satisfies condition (b).

We now remark that condition (b) is equivalent to condition (c) for a morphism
$f:A\rightarrow B$.
Indeed, the part of condition (c) for $i=-1$ is, by the definition of $\pi _0$,
identical to the condition  that $f$ induces a surjection $\pi _0(A)\rightarrow
\pi _0(B)$. And the remaining conditions for $i=0,\ldots , n+1$ are identical to
the  conditions of (c) corresponding to $j=i-1=-1,\ldots , (n-1)+1$ for all the
morphisms of $n-1$-groupoids  $Hom _A(x,y)\rightarrow Hom _B(f(x), f(y))$.
(In terms of $u$ and $v$ appearing in the condition in question, take $x$ to
be the source of the source of the source \ldots , and take $y$ to be the
target of the target of the target \ldots ).   Thus by induction on $n$
(i.e. by the equivalence $(b)\Leftrightarrow (c)$ for $n-1$-groupoids),
the conditions (c) for $f$ for $i=0,\ldots , n+1$, are equivalent to the
conditions that  $Hom _A(x,y)\rightarrow Hom _B(f(x), f(y))$
be equivalences of $n-1$-groupoids. Thus condition (c) for $f$
is equivalent to condition (b) for $f$, which completes the proof of part
(III) of the theorem.

\bigskip

We now proceed with the proof of part (I) of Theorem \ref{thmdef}. Note first of
all that the implications $(1)\Rightarrow (2)$ and $(2)\Rightarrow (3)$ are
easy. We give a short discussion of $(1)\Rightarrow (3)$ anyway, and then
we prove $(3)\Rightarrow (2)$ and $(2)\Rightarrow (1)$.

Note also that
the equivalence $(1)\Leftrightarrow (2)$ is the content of Proposition 1.6 of
\cite{KV}; we give a proof here because the proof of
Proposition 1.6 was ``left to the reader'' in \cite{KV}.

\medskip

\noindent
{\bf $(1)\Rightarrow (3)$:}\,\, Suppose $A$ is a strict $n$-category
satisfying condition $(1)$. This condition (from \cite{KV}) is compatible with
truncation, so $\tau _{\leq 1} (A)$ satisfies condition $(1)$ for
$1$-categories; which in turn is equivalent to the standard condition of being a
$1$-groupoid, so we get that $\tau _{\leq 1}(A)$ is a $1$-groupoid. On the other
hand, the conditions $(1)$ from \cite{KV} for $i$-arrows, $1\leq i \leq n$,
include the same conditions for the $i-1$-arrows  of $Hom _A(x,y)$ for any
$x,y\in Ob(A)$ (the reader has to verify this by looking at the definition in
\cite{KV}). Thus by the inductive statement of the present theorem for strict
$n-1$-categories, $Hom _A(x,y)$ is a strict $n-1$-groupoid. This shows that  $A$
satisfies condition $(3)$.

\medskip

\noindent
{\bf $(3)\Rightarrow (2)$:}\,\, Suppose $A$ is a strict $n$-category satisfying
condition $(3)$. It already satisfies the first part of condition $(2)$, by
hypothesis. Thus we have to show the second part, for example that for $f:
x\rightarrow y$ in $Ob(Hom _A(x,y))$, composition with $f$ induces an
equivalence
$$
Hom _A(y,z)\rightarrow Hom _A(x,z)
$$
(the other part is dual and has the same proof which we won't repeat here).

In order to prove this, we need to make a digression about the effect of
composition with $2$-morphisms. Suppose $f,g\in Ob (Hom _A(x,y))$ and
suppose that $u$ is a $2$-morphism from $f$ to $g$---this last supposition
may be
rewritten
$$
u\in Ob(Hom _{Hom _A(x,y)}(f,g)).
$$
{\em Claim:}
Suppose $z$ is another object; we claim that if composition with $f$
induces an equivalence $Hom _A(y,z)\rightarrow Hom _A(x,z)$, then composition
with $g$ also induces an equivalence $Hom _A(y,z)\rightarrow Hom _A(x,z)$.

To prove the claim, suppose that  $h,k$ are two $1$-morphisms from $y$ to $z$.
We now obtain a diagram
$$
\begin{array}{ccc}
Hom _{Hom _A(y,z)}(h,k) & \rightarrow & Hom _{Hom _A(x,z)}(hf, kf)\\
\downarrow && \downarrow \\
Hom _{Hom _A(x,z)}(hg, kg) & \rightarrow & Hom _{Hom _A(x,z)}(hf, kg),
\end{array}
$$
where the top arrow is given by composition $\ast _0$ with $1_f$; the left arrow
by composition $\ast _0$ with $1_g$; the bottom arrow by composition $\ast _1$
with the $2$-morphism $h\ast _0u$; and the right morphism is given by
composition
with $k\ast _0u$. This diagram commutes (that is the ``Godement rule'' or
``interchange rule'' cf \cite{KV} p. 32). By the inductive statement of the
present theorem (version (2) of the groupoid condition) for the $n-1$-groupoid
$Hom _A(x,z)$, the morphisms on the bottom and on the right in the above
diagram are equivalences. The hypothesis in the claim that $f$ is an
equivalence means that the morphism along the top of the diagram is an
equivalence; thus by the sub-lemma (part (IV) of the present theorem) applied
to the $n-2$-groupoids in the diagram, we get that the morphism on the left of
the diagram is an equivalence. This provides the second half of the criterion
(b) of part (III) for showing that the morphism of composition with $g$,
$Hom _A(y,z)\rightarrow Hom _A(x,z)$, is an equivalence of $n-1$-groupoids.

To finish the proof of the claim, we now verify the first half of criterion (b)
for the morphism of composition with $g$ (in this part we use directly the
condition (3) for $A$ and don't use either $f$ or $u$).  Note that $\tau _{\leq
1}(A)$ is a $1$-groupoid, by the condition (3) which we are assuming. Note also
that (by definition)
$$
\pi _0Hom _A(y,z)= Hom _{\tau _{\leq 1}A}(y,z) \;\;\; \mbox{and}\;\;\;
\pi _0Hom _A(x,z)= Hom _{\tau _{\leq 1}A}(x,z),
$$
and the morphism in question here is just the morphism of composition by the
image of $g$ in $\tau _{\leq 1}(A)$. Invertibility of this morphism in
$\tau _{\leq 1}(A)$ implies that the composition morphism
$$
Hom _{\tau _{\leq 1}A}(y,z)\rightarrow Hom _{\tau _{\leq 1}A}(x,z)
$$
is an isomorphism. This completes verification of the first half of criterion
(b), so we get that composition with
$g$ is an equivalence. This completes the proof of the claim.

We now return to the proof of the composition condition for (2). The fact that
$\tau _{\leq 1}(A)$ is a $1$-groupoid implies that given $f$ there is another
morphism $h$ from $y$ to $x$ such that the class  of $fh$ is equal to the
class of
$1_y$ in $\pi _0Hom _A(y,y)$, and the class of $hf$ is equal to
the class of $1_x$ in $\pi _0Hom _A(x,x)$. This means that there exist
$2$-morphisms $u$ from $1_y$ to $fh$, and $v$ from $1_x$ to $hf$. By the above
claim (and the fact that the compositions with $1_x$ and $1_y$ act as the
identity and in particular are equivalences), we get that composition
with $fh$ is an equivalence
$$
\{ fh\} \times Hom _A(y,z) \rightarrow Hom _A(y,z),
$$
and that composition with $hf$ is an equivalence
$$
\{ hf\} \times Hom _A(x,z)\rightarrow Hom _A(x,z).
$$
Let
$$
\psi _f: Hom _A(y,z)\rightarrow Hom _A(x,z)
$$
be the morphism of composition with $f$, and let
$$
\psi _h: Hom _A(x,z)\rightarrow Hom _A(y,z)
$$
be the morphism of composition with $h$. We have seen that $\psi _h\psi _f$ and
$\psi _f\psi _h$ are equivalences. By the second sub-lemma (part (V) of the
theorem) applied to $n-1$-groupoids, these imply that $\psi _f$ is an
equivalence.

The proof for composition in the other direction is the same; thus we have
obtained condition (2) for $A$.

\medskip

\noindent
{\bf $(2)\Rightarrow (1)$:}\,\, Look at the condition (1) by refering to
\cite{KV}: in question are the conditions $GR'_{i,k}$ and $GR''_{i,k}$ ($i<k\leq
n$) of Definition 1.1, p. 33 of \cite{KV}. By the inductive version of the
present equivalence for $n-1$-groupoids and by the part of condition (2) which
says that the $Hom _A(x,y)$ are $n-1$-groupoids, we obtain the conditions
$GR'_{i,k}$ and $GR''_{i,k}$ for $i\geq 1$. Thus we may now restrict our
attention to the condition $GR'_{0,k}$ and  $GR''_{0,k}$.  For a $1$-morphism
$a$ from $x$ to $y$, the conditions $GR'_{0,k}$ for all $k$ with respect to $a$,
are the same as the condition that for all $w$, the morphism of
pre-multiplication by $a$ $$
Hom _A(w,x)\times \{ a\} \rightarrow Hom _A(w,y)
$$
is an equivalence according to the version (c) of the notion of equivalence (cf
Part (III) of this theorem). Thus, condition $GR'_{0,k}$ follows from the
second part of condition (2) (for pre-multiplication). Similarly condition
$GR''_{0,k}$ follows from the second part of condition (2) for
post-multiplication by every $1$-morphism $a$. Thus condition (2) implies
condition (1). This completes the proof of Part (I) of the theorem.

\bigskip

For the sub-lemma (part (IV) of the theorem), using the fact that isomorphisms
of sets satisfy the same ``three for two'' property, and using the
characterization of equivalences in terms of homotopy groups (condition (a))
we immediately get two of the three statements: that if $f$ and $g$ are
equivalences then $gf$ is an equivalence; and that if $gf$ and $g$ are
equivalences then $f$ is an equivalence. Suppose now that $gf$ and $f$ are
equivalences; we would like to show that $g$ is an equivalence. First of all
it is clear that if $x\in Ob(A)$ then $g$ induces an isomorphism
$\pi _i(B, f(x))\cong \pi _0(C, gf(x))$ (resp. $\pi _0(B)\cong \pi _0(C)$).
Suppose now that $y\in Ob(B)$, and choose a $1$-morphism $u$ going from $y$ to
$f(x)$ for some $x\in Ob(A)$ (this is possible because $f$ is surjective on $\pi
_0$).  By condition (2) for being a groupoid, composition with
$u$ induces equivalences along the top
row of the diagram
$$
\begin{array}{ccccc}
Hom _B(y,y) &\rightarrow &Hom _B(y, f(x))&\leftarrow &Hom _B(f(x), f(x))\\
\downarrow && \downarrow && \downarrow \\
Hom _C(g(y),g(y)) &\rightarrow &Hom _C(g(y), gf(x))&\leftarrow &
Hom _C(gf(x), gf(x)).
\end{array}
$$
Similarly composition with $g(u)$ induces equivalences along the bottom row.
The sub-lemma for $n-1$-groupoids applied to the sequence
$$
Hom _A(x,x)\rightarrow Hom _B(f(x), f(x))\rightarrow Hom _C(gf(x), gf(x))
$$
as well as the hypothesis that $f$ is an equivalence, imply that the rightmost
vertical arrow in the above diagram is an equivalence. Again applying the
sub-lemma to these $n-1$-groupoids yields that the leftmost vertical arrow is
an equivalence. In particular $g$ induces isomorphisms
$$
\pi _i(B,y) = \pi _{i-1}(Hom _B(y,y), 1_y) \stackrel{\cong}{\rightarrow}
\pi _{i-1}(Hom _C(g(y), g(y)),1_{g(y)}) = \pi _i(C, g(y)).
$$
This completes the verification of condition (a) for the morphism $g$,
completing the proof of part (IV) of the theorem.

Finally we prove the second sub-lemma, part (V) of the theorem (from which we
now adopt the notations $A,B,C,D,f,g,h$). Note first of all that
applying $\pi _0$ gives the same situation for maps of sets, so $\pi _0(g)$ is
an isomorphism. Next, suppose $x\in Ob(A)$. Then we obtain a sequence
$$
\pi _i(A,x)\rightarrow \pi _i(B,f(x))
\rightarrow \pi _i(C, gf(x))\rightarrow \pi _i(D, hgf(x)),
$$
such that the composition of the first pair and also of the last pair are
isomorphisms; thus $g$ induces an isomorphism
$\pi _i(B , f(x))\cong \pi _i(C, gf(x))$. Now, by the same argument as for Part
(IV) above, (using the hypothesis that $f$ induces a surjection $\pi
_0(A)\rightarrow \pi _0(B)$) we get that for any object $y\in Ob(B)$, $g$
induces an isomorphism $\pi _i(B , y)\cong \pi _i(C, g(y))$. By definition (a)
of Part (III) we have now shown that $g$ is an equivalence.
This completes the proof of the theorem.
\eop

Let $nStrGpd$ be the category of strict $n$-groupoids.

We close out this section by looking at how the groupoid condition fits in with
the discussion of \ref{scholium} and \ref{iterate}. Let $C$ be a strict
$n$-category with only one object $x$. Then $C$ is an $n$-groupoid if and only
if $Hom _C(x,x)$ is an $n-1$-groupoid and
$\pi _0Hom _C(x,x)$ (which has a structure of monoid) is a group.
This is version (3) of the definition of groupoid in \ref{thmdef}. Iterating
this remark one more time we get the following statement.

\begin{lemma}
\label{scholiumgpd}
The construction of \ref{scholium} establishes an equivalence of categories
between the strict $n$-groupoids having only one object and only one
$1$-morphism, and the abelian monoid-objects $G$ in $(n-2)StrGpd$ such that the
monoid $\pi _0(G)$ is a group.
\end{lemma}
\eop

\begin{corollary}
\label{iterategpd}
Suppose $C$ is a strict $n$-category  having only one object and only one
$1$-morphism, and let $B$ be the strict $n+1$-category of \ref{iterate}
with one object $b$ and $Hom _B(b,b)= C$. Then $B$ is a strict $n+1$-groupoid if
and only if $C$ is a strict $n$-groupoid.
\end{corollary}
{\em Proof:}
Keep the notations of the proof of \ref{iterate}.
If $C$ is a groupoid this means that $G$ satisfies the condition that $\pi
_0(G)$ be a group, which in turn implies that $U$ is a groupoid. Note that $\pi
_0(U)=\ast$ is automatically a group; so applying the observation
\ref{scholiumgpd} once again, we get that $B$ is a groupoid. In the other
direction, if $B$ is a groupoid then $C=Hom _B(b,b)$ is a groupoid by versions
(2) and (3) of the definition of groupoid.
\eop

\numero{Realization functors}

Recall that $nStrGpd$ is the category of strict $n$-groupoids as defined above
\ref{thmdef}. Let $Top$ be the category of topological spaces. The following
definition encodes the minimum of what one would expect for a reasonable
realization functor from strict $n$-groupoids to spaces.

\begin{definition}
\label{realizationdef}
A {\em realization functor for strict $n$-groupoids} is a functor
$$
\Re : nStrGpd \rightarrow Top
$$
together with the following natural transformations:
$$
r:Ob (A) \rightarrow \Re (A);
$$
$$
\zeta _i(A,x): \pi _i (A, x) \rightarrow \pi _i (\Re (A), r(x)),
$$
the latter including
$\zeta _0(A): \pi _0(A)\rightarrow \pi _0(\Re (A))$; such that
the $\zeta _i(A,x)$ and $\zeta _0(A)$ are isomorphisms for $0\leq i \leq n$,
and such that the $\pi _i(\Re (A), y)$ vanish for $i>n$.
\end{definition}

\begin{theorem}
\label{realization}
{\rm (\cite{KV})}
There exists a realization functor $\Re$ for strict  $n$-groupoids.
\end{theorem}

Kapranov and Voevodsky \cite{KV} construct such a functor. Their construction
proceeds by first defining a notion of ``diagrammatic set''; they define a
realization functor from $n$-groupoids to diagrammatic sets (denoted $Nerv$),
and then define the topological realization of a diagrammatic set (denoted $|
\cdot |$). The composition of these two constructions gives a realization
functor  $$
G \mapsto \Re _{KV}(G):= | Nerv(G)|
$$
from
strict $n$-groupoids to spaces. Note that this functor $\Re_{KV}$ satisfies the
axioms of \ref{realizationdef} as a consequence of Propositions 2.7 and 3.5 of
\cite{KV}.

One obtains a different construction by considering strict $n$-groupoids as weak
$n$-groupoids in the sense of \cite{Tamsamani} (multisimplicial sets) and then
taking the realization of \cite{Tamsamani}.  This construction is actually
probably due to someone from the Australian school many years beforehand and we
call it the {\em standard realization} $\Re _{\rm std}$. The properties of
\ref{realizationdef} can be extracted from \cite{Tamsamani} (although again
they are probably classical results).

We don't claim here that any two realization functors must be the same, and
in particular the realization $\Re _{KV}$ could {\em a priori} be different from
the standard one. This is why we shall work, in what follows, with an arbitrary
realization functor satisfying the axioms of \ref{realizationdef}.

Here are some
consequences of the axioms for a realization functor. If $C\rightarrow C'$ is a
morphism of strict $n$-groupoids inducing isomorphisms on the $\pi _i$ then $\Re
(C)\rightarrow \Re (C')$ is a weak homotopy equivalence.
Conversely if $f:C\rightarrow C'$ is a morphism of strict $n$-groupoids which
induces a weak equivalence of realizations then $f$ was an equivalence.

\numero{The case of the standard realization}

Before getting to our main result which concerns an arbitrary realization
functor satisfying \ref{realizationdef}, we take note of an easier argument
which shows that the standard realization functor cannot give rise to arbitrary
homotopy types.

\begin{definition}
\label{compatiblelooping}
A collection of realization functors $\Re ^n$ for $n$-groupoids ($0\leq n <
\infty$) satisfying \ref{realizationdef} is said to be {\em compatible with
looping} if there exist transformations natural in an $n$-groupoid $A$ and an
object $x\in Ob(A)$,
$$
\varphi (A, x): \Re ^{n-1}(Hom _A(x,x))\rightarrow \Omega ^{r(x)}\Re ^n(A)
$$
(where $\Omega ^{r(x)}$ means the space of loops based at $r(x)$), such that
for $i\geq 1$ the following diagram commutes:
$$
\begin{array}{ccc}
\pi _i(A, x) & = \pi _{i-1}(Hom _A(x,x), 1_x) \rightarrow &
\pi _{i-1}(\Re ^{n-1}(Hom _A(x,x)), r(1_x))\\
\downarrow &&\downarrow \\
\pi _i(\Re ^n(A), r(x)) & \leftarrow & \pi _{i-1}(  \Omega ^{r(x)}\Re ^n(A),
cst(r(x)))
\end{array}
$$
where the top arrow is $\zeta _{i-1}(Hom _A(x,x), 1_x)$, the left arrow is
$\zeta _{i}(A,x)$, the right arrow is induced by $\varphi (A, x)$, and the
bottom arrow is the canonical arrow from topology. (When $i=1$, suppress the
basepoints in the $\pi _{i-1}$ in the diagram.)
\end{definition}

{\em Remark:} The arrows on the top, the bottom and the left are isomorphisms
in the above diagram, so the arrow on the right is an isomorphism and we obtain
as a corollary of the definition that the $\varphi (A,x)$ are actually weak
equivalences.

{\em Remark:} The collection of standard realizations $\Re ^n_{\rm std}$ for
$n$-groupoids, is compatible with looping. We leave this as an exercise for the
reader.

Recall the statements of \ref{iterate} and \ref{iterategpd}:
if $A$ is a strict  $n$-category with only one object $x$ and only one
$1$-morphism $1_x$, then there exists a strict $n+1$-category $B$ with
one object $y$, and with $Hom _B(y,y)=A$; and  $A$ is a strict  $n$-groupoid if
and only if $B$ is a strict $n+1$-groupoid.

\begin{corollary}
\label{forstandard}
Suppose $\{ \Re ^n\}$ is a collection of realization  functors
\ref{realizationdef} compatible with looping \ref{compatiblelooping}.
Then if $A$ is a $1$-connected strict $n$-groupoid (i.e. $\pi _0(A)=\ast$ and
$\pi _1(A,x)=\{ 1\}$), the space $\Re ^n(A)$ is weak-equivalent to a loop
space.
\end{corollary}
{\em Proof:}
Let $A'\subset A$ be the sub-$n$-category having one object $x$ and one
$1$-morphism $1_x$. For $i\geq 2$ the inclusion induces isomorphisms
$$
\pi _i(A', x) \cong \pi _i(A,x),
$$
and in view of the $1$-connectedness of $A$ this means (according to the
definition of \ref{thmdef} III (a)) that the morphism $A'\rightarrow A$ is an
equivalence. It follows (by definition \ref{realizationdef}) that $\Re
^n(A')\rightarrow \Re ^n(A)$ is a weak equivalence.  Now $A'$ satisfies the
hypothesis of \ref{iterate}, \ref{iterategpd} as recalled above, so there is an
$n+1$-groupoid $B$ having one object $y$ such that $A'=Hom _B(y,y)$. By the
definition of ``compatible with looping'' and the subsequent remark that the
morphism $\varphi (B,y)$ is a weak equivalence, we get that $\varphi (B,y)$
induces a weak equivalence
$$
\Re ^n(A') \rightarrow \Omega ^{r(y)}\Re ^{n+1}(B).
$$
Thus $\Re ^n(A)$ is weak-equivalent to the loop-space of $\Re ^{n+1}(B)$.
\eop

The following corollary is a statement which seems to be due to C. Berger
\cite{Berger} (although the statement appears without proof in Grothendieck
\cite{Grothendieck}). See also R. Brown and coauthors \cite{RBrown1}
\cite{BrownGilbert} \cite{BrownHiggins} \cite{BrownHiggins2}.

\begin{corollary}
\label{berger}
{\rm (C. Berger \cite{Berger})} There is no strict $3$-groupoid $A$ such that
the standard realization $\Re _{\rm std} (A)$ is weak-equivalent to the
$3$-type of $S^2$.
\end{corollary}
{\em Proof:}
The $3$-type of $S^2$ is not a loop-space. By the previous corollary
(and the fact that the standard realizations are compatible with looping, which
we have above left as an exercise for the reader), it is impossible for
$\Re _{\rm std}(A)$ to be the $3$-type of $S^2$.
\eop

\numero{Nonexistence of strict $3$-groupoids giving rise to the $3$-type
of $S^2$}

It is not completely clear whether Kapranov and Voevodsky claim that their
realization functors are compatible with looping in the sense of
\ref{compatiblelooping}, so Berger's negative result (Corollary \ref{berger}
above) might not apply. The main work of the present paper is to extend
this negative result to {\em any} realization functor satisfying the minimal
definition \ref{realizationdef}, in particular getting a result which applies
to the realization functor of \cite{KV}.

\begin{proposition}
\label{noS2}
Let $\Re$ be any realization functor satisfying the properties of Definition
\ref{realizationdef}. Then there does not exist a strict $3$-groupoid $C$ such
that $\Re (C)$ is weak-equivalent to the $3$-truncation of the homotopy type of
$S^2$.
\end{proposition}

\begin{corollary}
Let  $\Re _{KV}$ be the realization functor of Kapranov and Voevodsky
\cite{KV} cf the discussion above. If we assume that Propositions 2.7 and 3.5 of
\cite{KV} (stating that $\Re _{KV}$ satisfies the axioms \ref{realizationdef})
are true, then Corollary 3.8 of \cite{KV} is not true, i.e. $\Re_{KV}$ does not
induce an equivalence between the homotopy categories of strict
$3$-groupoids and
$3$-truncated topological spaces.
\end{corollary}
{\em Proof:}
According to Proposition \ref{noS2}, for any realization functor
satisfying \ref{realizationdef}, the induced functor on the homotopy categories
is not essentially surjective: its essential image doesn't contain the $3$-type
of $S^2$.
\eop

Proposition \ref{noS2} is very similar to the result of Brown and Higgins
\cite{BrownHiggins} and also the recent result of C. Berger \cite{Berger} (cf
\ref{berger} above). As was noted in \cite{KV}, the result of Brown and Higgins
concerns the more restrictive notion of groupoid where one requires that all
morphisms have strict inverses (however, see also \cite{RBrown1},
\cite{BrownHiggins2}). As in \cite{KV}, that restriction is not included in the
definition \ref{thmdef}. Berger considers strict $n$-groupoids according to the
definition \ref{thmdef} (i.e. with inverses non-strict) as well, but his
negative result applies only to a standard realization functor and as such,
doesn't {\em a priori} directly contradict \cite{KV}.

The basic difference in the present approach is that we
make no reference to any particular construction of $\Re$ but show that the
proposition holds for any realization construction having the properties of
Definition \ref{realizationdef}.

The fact that strict
$n$-groupoids don't model all homotopy types is also mentionned in Grothendieck
\cite{Grothendieck}. The basic idea in the setting of $3$-categories not
necessarily groupoids, is contained in some examples which G. Maltsiniotis
pointed out to me, in Gordon-Power-Street \cite{Gordon-Power-Street} where there
are given examples of weak $3$-categories not equivalent to strict ones.  This
in turn is related to the difference between braided monoidal categories and
symmetric monoidal categories, see for example the nice discussion in Baez-Dolan
\cite{BaezDolan}.

\bigskip

In order to prove \ref{noS2}, we will prove the following statement (which
contains the main part of the argument). It basically says that the Postnikov
tower of a simply connected strict $3$-groupoid $C$, splits.

\begin{proposition}
\label{diagramme}
Suppose $C$ is a strict $3$-groupoid with an object $c$ such that  $\pi
_0(C)=\ast$, $\pi _1(C,c)=\{ 1\}$, $\pi _2(C,c)
 = \zz$ and $\pi _3(C,c)=H$ for an abelian group $H$. Then there exists a
diagram of strict $3$-groupoids
$$
C \stackrel{g}{\leftarrow} B \stackrel{f}{\leftarrow} A
\stackrel{h}{\rightarrow} D
$$
with objects $b\in Ob(B)$, $a\in Ob(A)$, $d\in Ob(D)$ such that
$f(a)=b$, $g(b)=c$, $h(a)=d$. The diagram is
such that $g$ and $f$ are equivalences of strict $3$-groupoids, and such that
$\pi _0(D)=\ast$, $\pi _1(D,d)=\{ 1\}$, $\pi _2(D,d)=\{ 0\}$, and
such that $h$ induces an isomorphism
$$
\pi _3(h): \pi _3(A,a)=H \stackrel{\cong}{\rightarrow} \pi _3(D,d).
$$
\end{proposition}

\subnumero{Proof of Proposition \ref{noS2} using Proposition \ref{diagramme}}

Suppose for the moment that we know Proposition \ref{diagramme}; with this we
will prove \ref{noS2}. Fix a realization functor $\Re$ for strict
$3$-groupoids satisfying the axioms \ref{realizationdef}, and assume that $C$ is
a strict $3$-groupoid such that $\Re (C)$ is weak homotopy-equivalent to the
$3$-type of $S^2$. We shall derive a contradiction.

In {\em r\'esum\'e} the argument is this: that applying the realization functor
to the diagram given by  \ref{diagramme} and inverting the first two maps which
are weak homotopy equivalences, we would get a map
$$
\tau _{\leq 3}(S^2)= \Re (C) \rightarrow \Re (D) = K(H, 3)
$$
(with $H=\zz $).
This is a class in $H^3(S^2, H)$.
The hypothesis that $\Re (h)$ is an isomorphism on $\pi _3$ means that this
class is nonzero when applied to $\pi _3(S^2)$ via the Hurewicz homomorphism;
but $H^3(S^2, \zz )= 0$, a contradiction.

Here is a full description of the argument.
Apply Proposition \ref{diagramme} to $C$. Choose an object $c\in Ob(C)$. Note
that, because of the isomorphisms between homotopy sets or groups
\ref{realizationdef}, we have $\pi _0(C)=\ast$,
$\pi _1(C,c)=\{ 1\}$, $\pi _2(C,c)
 = \zz$ and $\pi _3(C,c)=\zz $, so \ref{diagramme} applies with $H=\zz$.
We obtain a sequence of strict $3$-groupoids
$$
C \stackrel{g}{\leftarrow} B \stackrel{f}{\leftarrow} A
\stackrel{h}{\rightarrow} D.
$$
This gives the diagram of spaces
$$
\Re (C) \stackrel{\Re (g)}{\leftarrow} \Re (B) \stackrel{\Re
(f)}{\leftarrow} \Re
(A)  \stackrel{\Re (h)}{\rightarrow} \Re (D).
$$
The axioms \ref{realizationdef} for $\Re$ imply that $\Re$ transforms
equivalences of strict $3$-groupoids into weak homotopy equivalences of spaces.
Thus $\Re (f)$ and $\Re (g)$ are weak homotopy equivalences and we get that
$\Re (A)$ is weak homotopy equivalent to the $3$-type of $S^2$.

On the other hand, again by the
axioms \ref{realizationdef}, we have that $\Re (D)$ is $2$-connected, and $\pi
_3(\Re (D), r(d))=H$ (via the isomorphism $\pi _3(D,d)\cong H$ induced by $h$,
$f$ and $g$).
By the Hurewicz theorem there is a class $\eta \in H^3(\Re (D),
H)$ which induces an isomorphism
$$
{\bf Hur}(\eta ): \pi _3(\Re (D), r(d))\stackrel{\cong}{\rightarrow} H.
$$
Here
$$
{\bf Hur} : H^3(X , H)\rightarrow Hom (\pi _3(X,x), H)
$$
is the Hurewicz map for any pointed space $(X,x)$; and the cohomology is
singular cohomology (in particular it only depends on the weak homotopy type of
the space).

Now look at the pullback of this class
$$
\Re (h)^{\ast}(\eta )\in H^3(\Re (A), H).
$$
The hypothesis that $\Re (u)$ induces an isomorphism on $\pi _3$ implies that
$$
{\bf Hur}(\Re (h)^{\ast}(\eta )): \pi _3(\Re (A),
r(a))\stackrel{\cong}{\rightarrow} H.
$$
In particular, ${\bf Hur}(\Re (h)^{\ast}(\eta ))$ is nonzero so
$\Re (h)^{\ast}(\eta )$ is nonzero in $H^3(\Re (A), H)$. This is a
contradiction because $\Re (A)$ is weak homotopy-equivalent to the $3$-type
of $S^2$, and $H=\zz$, but $H^3(S^2 , \zz )=\{ 0 \}$.

This contradiction completes the proof of Proposition \ref{noS2}, assuming
Proposition \ref{diagramme}.
\eop

\subnumero{Proof of Proposition \ref{diagramme}}

This is the main part of the argument. We start with a strict groupoid $C$ and
object $c$, satisfying the hypotheses of \ref{diagramme}.

The first step is to construct $(B,b)$. We let $B\subset C$ be the
sub-$3$-category having only one object $b=c$, and only one $1$-morphism
$1_b=1_c$. We set
$$
Hom _{Hom _B(b,b)}(1_b, 1_b):=Hom _{Hom _C(c,c)}(1_c, 1_c) ,
$$
with the same composition law.
The map $g: B\rightarrow C$ is the inclusion.

Note first of all that $B$ is a strict
$3$-groupoid. This is easily seen using version (1) of the definition
\ref{thmdef} (but one has to look at the conditions in \cite{KV}).
We can also verify it using condition (3). Of course $\tau _{\leq 1}(B)$ is the
$1$-category with only one object and only one morphism, so it is a groupoid.
We have to verify that $Hom _B(b,b)$ is a strict $2$-groupoid. For this, we
again apply condition (3) of \ref{thmdef}. Here we note that
$$
Hom _B(b,b)\subset Hom _C(c,c)
$$
is the full sub-$2$-category with only one object $1_b=1_c$. Therefore, in view
of the definition of $\tau _{\leq 1}$, we have that
$$
\tau _{\leq 1}Hom _B(b,b)\subset \tau _{\leq 1}Hom _C(c,c)
$$
is a full subcategory. A full subcategory of a $1$-groupoid is again a
$1$-groupoid, so $\tau _{\leq 1}Hom _B(b,b)$ is a $1$-groupoid. Finally,
$Hom _{Hom _B(b,b)}(1_b, 1_b)$ is a $1$-groupoid since by construction it
is the same as  $Hom _{Hom _C(c,c)}(1_c, 1_c)$ (which is a groupoid by condition
(3) applied to the strict $2$-groupoid $Hom _C(c,c)$). This shows that
$Hom _B(b,b)$ is a strict $2$-groupoid an hence that $B$ is a strict
$3$-groupoid.

Next, note that $\pi _0(B)=\ast$ and $\pi _1(B,b)=\{ 1\}$. On the other hand,
for $i=2,3$ we have
$$
\pi _i(B,b)= \pi _{i-2}(Hom _{Hom _B(b,b)}(1_b, 1_b), 1^2_b)
$$
and similarly
$$
\pi _i(C,c)= \pi _{i-2}(Hom _{Hom _C(c,c)}(1_c, 1_c), 1^2_c),
$$
so the inclusion $g$ induces an equality $\pi _i(B,b) \stackrel{=}{\rightarrow}
\pi _i(C,c)$. Therefore, by definition (a) of equivalence \ref{thmdef}, $g$ is
an equivalence of strict $3$-groupoids. This completes the construction and
verification for $B$ and $g$.

Before getting to the construction of $A$ and $f$, we analyze the strict
$3$-groupoid $B$ in terms of the discussion of \ref{scholium} and
\ref{scholiumgpd}. Let
$$
G:= Hom _{Hom _B(b,b)}(1_b, 1_b).
$$
It is an abelian monoid-object in the category of $1$-groupoids, with abelian
operation denoted by $+: G\times G\rightarrow G$ and unit element denoted $0\in
G$ which is the same as $1_b$. The operation $+$ corresponds to both of the
compositions $\ast _0$ and $\ast _1$ in $B$.

The hypotheses on the homotopy
groups of $C$ also hold for $B$ (since $g$ was an equivalence). These translate
to the statements that $(\pi _0(G), +) = \zz$ and $Hom _G(0,0)=H$.

We now construct $A$ and $f$ via \ref{scholium} and
\ref{scholiumgpd}, by constructing a morphism $(G',+)\rightarrow (G,+)$ of
abelian monoid-objects in the category of $1$-groupoids. We do this by a type
of ``base-change'' on the monoid of objects, i.e. we will first define a
morphism $Ob(G')\rightarrow Ob(G)$ and then define $G'$ to be the groupoid
with object set $Ob(G')$ but with morphisms corresponding to those of $G$.

To accomplish the ``base-change'', start with the following construction. If
$S$ is a set, let ${\bf E}(S)$ denote the groupoid with $S$ as set of
objects, and with exactly one morphism between each pair of objects. If $S$ has
an abelian monoid structure then ${\bf E}(S)$ is an abelian monoid object in
the category of groupoids.

Note that for any groupoid $U$ there is a morphism of groupoids
$$
U\rightarrow {\bf E}(Ob(U)),
$$
and by ``base change'' we mean the following operation: take a set $S$ with a
map $p:S\rightarrow Ob(U)$ and look at
$$
V:= {\bf E}(S)\times _{{\bf E}(Ob(U))}U.
$$
This is a groupoid with $S$ as set of objects, and with
$$
Hom _V(s,t)= Hom _U(p(s), p(t)).
$$
If $U$ is an abelian monoid object in the category of groupoids, if $S$ is an
abelian monoid and if $p$ is a map of monoids then $V$ is again an abelian
monoid object in the category of groupoids.

Apply this as follows. Starting with $(G,+)$ corresponding to $B$
via \ref{scholium} and
\ref{scholiumgpd} as above,
choose objects $a,b \in Ob(G)$ such that the image of $a$ in $\pi
_0(G)\cong \zz$
corresponds to $1\in \zz$, and such that the image of $b$ in $\pi _0(G)$
corresponds to $-1\in \zz$. Let $N$ denote the abelian monoid, product of two
copies of the natural numbers, with objects denoted $(m,n)$ for nonnegative
integers $m,n$. Define a map  of abelian monoids
$$
p:N \rightarrow Ob(G)
$$
by
$$
p(m,n):= m\cdot a + n\cdot b := a+a+\ldots +a \, + \, b+b+\ldots +b.
$$
Note that this induces the surjection $N\rightarrow \pi _0(G)=\zz$
given by $(m,n)\mapsto m-n$.

Define $(G',+)$ as the base-change
$$
G':= {\bf E}(N) \times _{{\bf E}(Ob(G))} G,
$$
with its induced abelian monoid operation $+$. We have
$$
Ob (G')= N,
$$
and the second projection $p_2: G'\rightarrow G$ (which induces $p$ on object
sets) is fully faithful i.e.
$$
Hom _{G'}((m,n), (m',n'))= Hom _G(p(m,n), p(m',n')).
$$
Note that $\pi _0(G')=\zz$ via the map induced by $p$ or equivalently $p_2$.
To prove this, say that: (i) $N$ surjects onto $\zz$ so the map induced by
$p$ is
surjective; and (ii) the fact that $p_2$ is fully faithful implies that the
induced map $\pi _0(G')\rightarrow \pi _0(G)=\zz$ is injective.

We let $A$ be the strict $3$-groupoid corresponding to $(G',+)$
via \ref{scholium}, and let $f: A\rightarrow B$ be the
map corresponding to $p_2: G'\rightarrow G$ again via \ref{scholium}.
Let $a$ be the unique object of $A$ (it is mapped by $f$ to the unique object
$b\in Ob(B)$).

The fact that $(\pi _0(G'),+)=\zz$ is a group implies that $A$ is a
strict $3$-groupoid (\ref{scholiumgpd}). We have $\pi _0(A)=\ast$ and $\pi
_1(A,a)=\{ 1\}$. Also,
$$
\pi _2(A,a)= (\pi _0(G'), +) = \zz
$$
and $f$ induces an isomorphism from here to $\pi _2(B,b)=(\pi _0(G), +)=\zz$.
Finally (using the notation $(0,0)$ for the unit object of $(N,+)$ and the
notation $0$ for the unit object of $Ob(G)$),
$$
\pi _3(A,a)= Hom _{G'}((0,0),(0,0)),
$$
and similarly
$$
\pi _3(B,b)=Hom _G(0,0)=H;
$$
the map $\pi _3(f): \pi _3(A,a)\rightarrow \pi _3(B,b)$
is an isomorphism because it is the same as the map
$$
Hom _{G'}((0,0),(0,0))\rightarrow Hom _G(0,0)
$$
induced by $p_2: G'\rightarrow G$, and $p_2$ is fully faithful.
We have now completed the verification that $f$ induces isomorphisms on the
homotopy groups, so by  version (a) of the definition of equivalence
\ref{thmdef}, $f$ is an equivalence of strict $3$-groupoids.

We now construct $D$ and define the map $h$ by an explicit calculation in
$(G',+)$. First of all, let $[H]$ denote the $1$-groupoid with one object
denoted $0$, and with $H$ as group of endomorphisms:
$$
Hom _{[H]}(0,0):= H.
$$
This has a structure of abelian monoid-object in the category
of groupoids,  denoted $([H], +)$, because $H$ is an abelian group.
Let $D$ be the strict $3$-groupoid corresponding to $([H], +)$ via
\ref{scholium} and \ref{scholiumgpd}. We will construct a morphism
$h: A\rightarrow D$ via \ref{scholium} by constructing a morphism of abelian
monoid objects in the category of groupoids,
$$
h:(G', +)\rightarrow ([H], +).
$$
We will construct this morphism so that it induces the identity  morphism
$$
Hom _{G'}((0,0), (0,0))=H \rightarrow Hom _{[H]}(0,0)=H.
$$
This will insure that the morphism $h$ has the property required for
\ref{diagramme}.

The object $(1,1)\in N$ goes to $0\in \pi _0(G')\cong \zz$. Thus we may choose
an isomorphism $\varphi : (0,0)\cong (1,1)$ in $G'$.  For any $k$ let $k\varphi$
denote the isomorphism $\varphi + \ldots +\varphi$ ($k$ times) going from
$(0,0)$ to $(k,k)$.  On the other hand, $H$ is the automorphism group of
$(0,0)$ in $G'$. The operations $+$ and composition coincide on $H$. Finally,
for any $(m,n)\in N$ let $1_{m,n}$ denote the identity automorphism of the
object $(m,n)$.  Then any arrow $\alpha$ in $G$ may be uniquely written in the
form
$$
\alpha = 1_{m,n} + k\varphi + u
$$
with $(m,n)$ the source of $\alpha$, the target being $(m+k, n+k)$, and where
$u\in H$.

We have the following formulae for the composition $\circ$ of arrows in $G'$.
They all come from the basic rule
$$
(\alpha \circ \beta ) + (\alpha ' \circ \beta ')=
(\alpha + \alpha ') \circ (\beta + \beta ')
$$
which in turn comes simply from the fact that $+$ is a morphism of groupoids
$G'\times G'\rightarrow G'$ defined on the cartesian product of two copies of
$G$. Note in a similar vein that $1_{0,0}$ acts as the identity for the
operation
$+$ on arrows, and also that
$$
1_{m,n} + 1_{m',n'} = 1_{m+m', n+n'}.
$$

Our first equation is
$$
(1_{l,l} +k\varphi )\circ l\varphi = (k+l)\varphi .
$$
To prove this note that $l\varphi + 1_{0,0}= l\varphi$ and our basic formula
says
$$
(1_{l,l}\circ l_{\varphi} ) + (k\varphi \circ 1_{0,0})
=
(1_{l,l} +k\varphi )\circ (l\varphi + 1_{0,0} )
$$
but the left side is just $l\varphi + k\varphi = (k+l)\varphi$.

Now our basic formula, for a composition starting with
$(m,n)$, going first to $(m+l,n+l)$, then going to $(m+l+k, n+l+k)$, gives
$$
(1_{m+l,n+l} + k\varphi + u)\circ (1_{m,n} + l\varphi + v)
$$
$$
= (1_{m,n} + 1_{l,l} + k\varphi + u)\circ (1_{m,n} + l\varphi + v)
$$
$$
= 1_{m,n}\circ 1_{m,n}  + (1_{l,l} +k\varphi )\circ l\varphi
+ u\circ v
$$
$$
= 1_{m,n} + (k+l)\varphi + (u\circ v)
$$
where of course $u\circ v=u+v$.

This formula shows that the morphism $h$ from arrows of $G'$ to the group $H$,
defined by
$$
h(1_{m,n} + k\varphi + u):= u
$$
is compatible with composition.  This implies that it provides a morphism of
groupoids $h:G\rightarrow [H]$ (recall from above that $[H]$ is defined to
be the
groupoid with one object whose automorphism group is $H$).  Furthermore the
morphism $h$ is obviously compatible with the operation $+$ since
$$
(1_{m,n} + k\varphi + u)+ (1_{m',n'} + k'\varphi + u')=
$$
$$
(1_{m+m',n+n'} + (k+k')\varphi + (u+u'))
$$
and once again $u+u'=u\circ u'$ (the operation $+$ on $[H]$ being given by the
commutative operation $\circ$ on $H$).

This completes the construction of a morphism
$h: (G, +)\rightarrow ([H], +)$ which induces the identity on $Hom (0,0)$.
This corresponds to a morphism of strict $3$-groupoids $h: A\rightarrow D$
as required to complete the proof of Proposition \ref{diagramme}.
\eop

\bigskip

\pagebreak[4]

\numero{A remark on strict $\infty$-groupoids}

The nonexistence result of \ref{noS2} holds also for
strict $\infty$-groupoids as
defined in \cite{KV}.
Recall that Kapranov-Voevodsky \cite{KV} extend the notion of strict
$n$-category and strict $n$-groupoid to the case $n=\infty$. The definition is
made using condition (1), and the notion of equivalence is defined using (a) in
\ref{thmdef}. Note that the other characterizations of \ref{thmdef} don't
actually make sense in the case $n=\infty$ because they are inductive on $n$.

The only thing we need to know about the case $n=\infty$ is that there are
homotopy groups $\pi _i(A,a)$ of a strict $\infty$-groupoid $A$, and there are
truncation operations on strict $\infty$-groupoids such that $\tau _{\leq n}(A)$
is a strict $n$-groupoid  with a natural morphism
$$
A\rightarrow \tau _{\leq n}(A)
$$
inducing isomorphisms on homotopy groups for $i\leq n$.
(Here the $n$-groupoid $\tau _{\leq n}(A)$ is considered as an
$\infty$-groupoid in the obvious way.) The homotopy groups and truncation are
defined as in \cite{KV}---again, one has to avoid those versions of the
definitions \ref{thmdef} which are recursive on $n$.

We can extend the definition of \ref{realizationdef} to the case $n=\infty$.
It is immediate that for any realization functor $\Re$
satisfying the axioms \ref{realizationdef} for $n=\infty$, the morphism
$$
\Re (A)\rightarrow \Re (\tau _{\leq n}A)
$$
is the Postnikov truncation of $\Re (A)$. Applying \ref{noS2}, we obtain the
following result.

\begin{corollary}
\label{noInfiniteS2}
For any realization functor $\Re$ satisfying the axioms \ref{realizationdef} for
$n=\infty$, there does not exist a strict  $\infty$-groupoid $A$ (as defined by
Kapranov-Voevodsky \cite{KV}) such that $\Re (A)$ is weak
homotopy-equivalent to the $2$-sphere $S^2$.
\end{corollary}
{\em Proof:}
Note that if $\Re$ is a realization functor satisfying \ref{realizationdef}
for $n=\infty$, then composing with the inclusion $i_3^{\infty}$ from
the category of strict $3$-groupoids to the category of strict $\infty$-groupoids
we obtain a realization functor $\Re i_3^{\infty}$ for strict $3$-groupoids,
again satisfying \ref{realizationdef}.  If $A$ is a strict
$\infty$-groupoid then
the above truncation morphism, written more precisely, is
$$
A\rightarrow i_3^{\infty} \tau _{\leq 3}(A).
$$
This induces isomorphisms on the $\pi _i$ for $i\leq 3$. Applying $\Re$ we get
$$
\Re (A) \rightarrow \Re i_3^{\infty} \tau _{\leq 3}(A),
$$
inducing an isomorphism on homotopy groups for $i\leq 3$. In particular, if
$\Re (A)$ were weak homotopy-equivalent to $S^2$ then  this would imply that
$\Re i_3^{\infty} \tau _{\leq 3}(A)$ is the $3$-type of $S^2$.
In view of the fact that $\Re i_3^{\infty}$ is a realization functor according
to \ref{realizationdef} for strict $3$-groupoids, this would contradict
\ref{noS2}. Thus we conclude that there is no strict $\infty$-groupoid $A$ with
$\Re (A)$ weak homotopy-equivalent to $S^2$.
\eop

\numero{Conclusion}

One really needs to look at
some type of weak $3$-categories in order to get a hold of $3$-truncated
homotopy
types. O. Leroy \cite{Leroy} and apparently, independantly, Joyal and
Tierney \cite{JoyalTierney} were the first to do this. See also
Gordon, Power, Street \cite{Gordon-Power-Street} and Berger \cite{Berger} for
weak $3$-categories and $3$-types. Baues \cite{Baues} showed that
$3$-types correspond to {\em quadratic modules} (a generalization of
the notion of crossed complex) \cite{Baues}. Tamsamani \cite{Tamsamani} was the
first to relate weak $n$-groupoids and homotopy $n$-types. For other notions  of
weak $n$-category, see \cite{BaezDolanLetter} \cite{BaezDolanIII}
\cite{Batanin},
\cite{Batanin2}.

From homotopy theory (cf \cite{Lewis}) the following
type of yoga seems to come out: that it suffices to weaken any one of the
principal structures involved.  Most weak notions of $n$-category involve a
weakening of the associativity, or eventually of the Godement (commutativity)
conditions.

It seems likely that the arguments of \cite{KV} would show that
one could instead weaken the condition of being {\em unary} (i.e. having
identities for the operations) and keep associativity and Godement.  We give
a proposed definition of what this would mean and then state two conjectures.

\subnumero{Motivation}

Before giving the definition,
we motivate these remarks by looking at the {\em
Moore loop space} $\Omega ^x_M(X)$ of a space $X$ based at $x\in X$ (the Moore
loop space is referred to in \cite{KV} as a motivation for their
construction). Recall that $\Omega ^x_M(X)$ is the space of {\em pairs} $(r,
\gamma )$ where  $r$ is a real number $r\geq 0$ and $\gamma = [0,r]\rightarrow
X$ is a path starting and ending at $x$. This has the advantage of being a
strictly associative monoid. On the other side of the coin, the ``length''
function
$$
\ell : \Omega ^x_M(X)\rightarrow [0,\infty )\subset \rr
$$
has a special behavoir over $r=0$. Note that over the open half-line
$(0,\infty )$ the length function $\ell$ is a fibration (even a fiber-space)
with fiber homeomorphic to the usual loop space. However, the fiber over $r=0$
consists of a single  point, the constant path $[0,0]\rightarrow X$ based at
$x$. This additional point (which is the unit element of the monoid
$\Omega ^x_M(X)$) doesn't affect the topology of $\Omega ^x_M$ (at least if $X$
is locally contractible at $x$) because it is glued in as a limit of paths which
are more and more concentrated in a neighborhood of $x$. However, the map
$\ell$ is no longer a fibration over a neighborhood of $r=0$.  This is a bit of
a problem because $\Omega ^x_M$ is not compatible with direct products of the
space $X$; in order to obtain a compatibility one has to take the fiber product
over $\rr$ via the length function:
$$
\Omega ^{(x,y)}_M(X\times Y)= \Omega ^x_M(X) \times _{\rr} \Omega ^y_M(Y),
$$
and the fact that $\ell$ is not a fibration could end up causing a problem in
an attempt to iteratively apply a construction like the Moore loop-space.

Things seem to get better if we restrict to
$$
\Omega ^x_{M'}(X):=\ell ^{-1}((0,\infty ))\subset \Omega ^x_M(X) ,
$$
but this associative  monoid no longer has a strict unit. Even so, the constant
path of any positive length gives a weak unit.

A motivation coming from a different direction was an observation made by Z.
Tamsamani early in the course of doing his thesis. He was trying to define a
strict $3$-category $2Cat$ whose objects would be the strict $2$-categories and
whose morphisms would be the weak $2$-functors between $2$-categories (plus
notions of weak natural transformations and $2$-natural transformations).
At some point he came to the conclusion that one could adequately define
$2Cat$ as a strict $3$-category except that he couldn't get strict identities.
Because of this problem we abandonned the idea and looked toward weakly
associative $n$-categories. In retrospect it would be interesting to pursue
Tamsamani's construction of a strict $2Cat$ but with only weak identities.

\subnumero{Snucategories}

Now we get back to looking at what it could mean to weaken the unit property
for strict $n$-categories or strict $n$-groupoids.
We will define a notion of {\em $n$-snucategory}
(the initial `s' stands for strict, `nu' stands for non-unary) by induction on
$n$. There will be a notion of direct product of $n$-snucategories.
Suppose we know what these mean for $n-1$. Then an $n$-snucategory $C$ consists
of a set $C_0$ of objects together with, for every pair of objects $x,y\in C_0$
an $n-1$-snucategory $Hom _C(x,y)$ and composition morphisms
$$
Hom _C(x,y)\times Hom _C(y,z) \rightarrow Hom _C(x,z)
$$
which are strictly associative, such that the {\em weak unary condition} is
satisfied.  We now explain this condition. An element $e_x\in Hom _C(x,x)$ is
called a weak identity if:
\newline
---composition with $e$ induces  equivalences of $n-1$-snucategories
$$
Hom _C(x,y)\rightarrow Hom_C(x,y) , \;\;\;
Hom _C(y,x)\rightarrow Hom_C(y,x);
$$
---and if $e\cdot e$ is equivalent to $e$.

In order to complete the recursive definition we must define the notion of when
a morphism of $n$-snucategories is an equivalence, and we must define what it
means for two objects to be equivalent. A morphism is said to be an equivalence
if the induced morphisms on $Hom$ are equivalences of $n-1$-snucategories and
if it is essentially surjective on objects: each object in the target is
equivalent to the image of an object. It thus remains just to be seen what
equivalence of objects means. For this we introduce the {\em truncations}
$\tau _{\leq i}C$ of an $n$-snucategory $C$. Again this is done in the same way
as usual: $\tau _{\leq i}C$ is the $i$-snucategory with the same objects as $C$
and whose $Hom$'s are the truncations
$$
Hom_{\tau _{\leq i}C}(x,y):=\tau _{\leq i-1}Hom _C(x,y).
$$
This works for $i\geq 1$ by recurrence, and for $i=0$ we define the truncation
to be the set of isomorphism classes in $\tau _{\leq 1}C$.  Note that
truncation is compatible with direct product (direct products are defined in
the obvious way) and takes equivalences to equivalences. These statements used
recursively allow us to show that the truncations themselves satisfy the weak
unary condition. Finally, we say that two objects are equivalent if they map
to the same thing in $\tau _{\leq 0}C$.

Proceeding in the same way as in \S 2 above, we can define the
notion of $n$-snugroupoid.

\begin{conjecture}
There are functors $\Pi _n$ and $\Re$ between the categories of
$n$-snugroupoids and $n$-truncated spaces (going in the usual directions)
together with adjunction morphisms inducing an equivalence between the
localization of $n$-snugroupoids by equivalences, and $n$-truncated spaces by
weak equivalences.
\end{conjecture}

I think that the argument of \cite{KV} (which is unclear on
the question of identity elements) actually serves to prove the above
statement. I have called the above statement a ``conjecture'' because I
haven't checked this.

One might go out on a limb a bit more and make the following

\begin{conjecture}
The localization of the category of $n$-snucategories by equivalences
is equivalent to the localizations of the categories of weak $n$-categories
of Tamsamani and/or Baez-Dolan and/or Batanin by equivalences.
\end{conjecture}

This of course is of a considerably more speculative nature.

{\bf Caveat}: the above definition of ``snucategory'' is invented in an
{\em ad hoc} way, and in particular one naturally wonders whether or not the
equivalences $e\cdot e \sim e$ and higher homotopical data going along with
that, would need to be specified in order to get a good definition. I have no
opinion about this (the above definition being just the easiest thing to say
which gives some idea of what needs to be done).  Thus it is not completely
clear that the above definition of $n$-snucategory is the ``right'' one to fit
into the conjectures.

\end{document}